\documentclass[sn-mathphys,Numbered]{sn-jnl}

\usepackage{graphicx}%
\usepackage{multirow}%
\usepackage{amsmath,amssymb,amsfonts}%
\usepackage{amsthm}%
\usepackage{mathrsfs}%
\usepackage[title]{appendix}%
\usepackage{xcolor}%
\usepackage{textcomp}%
\usepackage{manyfoot}%
\usepackage{booktabs}%
\usepackage{algorithm}%
\usepackage{algorithmicx}%
\usepackage{algpseudocode}%
\usepackage{listings}%
\usepackage{enumerate}
\usepackage{setspace}

\theoremstyle{thmstyleone}%
\newtheorem{theorem}{Theorem}[section]%
\newtheorem{lemma}{Lemma}[section]%

\theoremstyle{thmstyletwo}%
\newtheorem{remark}{Remark}[section]%

\theoremstyle{thmstylethree}%
\newtheorem{definition}{Definition}[section]%

\theoremstyle{thmstylefour}%
\numberwithin{equation}{section}
\begin{document}
	
	\title[Heat flow gradient estimates under Finsler $CD(-K,N)$ geometric flow]{Gradient estimates and Harnack inequality for positive solutions to the heat equation under compact Finsler $CD(-K,N)$ geometric flows}
	
	
	\author*[1]{\fnm{Bin} \sur{Shen}}\email{shenbin@seu.edu.cn}
	
	
	
	\affil*[1]{\orgdiv{School of Mathematics}, \orgname{Southeast University}, \orgaddress{\street{Dongnandaxue Road 2}, \city{Nanjing}, \postcode{211189}, \state{Jiangsu}, \country{China}}}


	
	\abstract{In this manuscript, we extend the global gradient estimates for positive solutions to the heat equation under a general compact Finsler $CD(-K,N)$ geometric flow and derive the corresponding Harnack inequality. 
	}

	\keywords{Finsler geometric flow, Finsler Ricci flow, gradient estimates, heat flow, $CD(-K,N)$ condition, Harnack inequality}
	
	
	\pacs[MSC Classification]{58J35, 53C60, 53E20, 35R01}
	
	\maketitle

	\section{Introduction}
	
	Gradient estimates are instrumental in the study of elliptic and parabolic equations. They played a pivotal role in Perelman's groundbreaking work on the Poincar\'e conjecture and the geometrization conjecture, where he introduced a nonlinear heat equation for scalar curvature \cite{Perelman:2002}. L. Ni demonstrated the existence of Hermitian-Einstein metrics on Hermitian vector bundles over various complete noncompact K\"ahler manifolds by reducing the problem to finding nonnegative solutions to related Poisson equations \cite{Ni:2002}. Subsequently, Ni and Tam investigated the common solutions of the Poincar\'e-Lelong equation and the heat equation to analyze the K\"ahler-Ricci flow, providing several applications, including the long-time existence of solutions \cite{Ni-Tam:2004}. These methods can also be applied to obtain geometric and analytic bounds for the generalized Ricci flow \cite{Li:2024}. On the other hand, the study of geometric analysis on general $CD(K,N)$ spaces benefits from the rapid development of nonlinear analysis on Finsler metric measure spaces \cite{Kristaly:2023}.
	
	The Finsler heat equation on a Finsler metric measure space $(M,F,\mu)$ is 
	\begin{equation}\label{2}
		\partial_t u=\Delta u,
	\end{equation}
	where $\Delta=\Delta_{\mu}$ is a nonlinear Laplacian related to the metric and the measure.
	S. Ohta and K.-T. Sturm first established Li-Yau gradient estimates for positive solutions to \eqref{2} on compact Finsler metric measure spaces \cite{Ohta:2014}. This result was subsequently extended by Q. Xia to forward complete spaces \cite{Xia:2023} by using the iteration method intruduced in \cite{XiaC:2014}. In a separate development, D. Bao introduced and investigated the Finsler-Ricci flow equation through the framework of Akbar-Zadeh's Ricci tensor \cite{Bao:2007}, which is expressed as
    $$\frac{\partial}{\partial t}g=-2\mathcal{R}ic,$$
    where $\mathcal{R}ic(x,y)=Ric_{ij}(x,y)dx^i\otimes dx^j$, with $Ric_{ij}(x,y)=\left(\frac12\mathrm{Ric}(y)\right)_{y^iy^j}$. 	
	S. Lakzian investigated differential Harnack estimates for positive solutions to the heat equation under the Finsler-Ricci flow $(M,F(t),\mu)$ under the assumption of vanishing S-curvature \cite{Lakzian:2015}. This work was further re-studied by F. Zeng and Q. He in their series of studies \cite{Zeng-He:2019,Zeng:2020,Zeng:2021}, maintaining the same curvature assumption. Besides, X. Cheng revisited this problem in \cite{Cheng:2022} and removed the S-curvature vanishing assumption. In a broader context, S. Azami made an attempt by extending these results to a class of nonlinear parabolic equations under compact Finsler geometric flow \cite{Azami:2023}. However, the S-curvature conditions imposed in these works are unnecessarily strong (cf. Remark \ref{rmk2}), and more fatally, to the best of our knowledge, all existing literature on this problem in the Finsler setting contains critical errors (cf. \cite{Lakzian:2015,Zeng-He:2019,Zeng:2020,Zeng:2021,Cheng:2022,Azami:2023} and etc.). It is crucial to emphasize that all these results fundamentally  rely on the assumption that the distortion $\tau$ is Chern horizontally invariant (cf. Remark \ref{rmk1}), a condition that significantly constrains their applicability. Given these issues, a rigorous correction of these estimates is essential for the proper development of geometric analysis on metric measure spaces under geometric flow.
	
	In this manuscript, to prove the differential Harnack inequalities under general Finsler geometric flow, we utilize the remarkable $CD(-K,N)$ condition to overcome the obstacles caused by the S-curvature and the measure $\mu$. 
	
	Let ($M$,$F$,$\mu$) be a compact Finsler metric measure space. The Finsler geometric flow is 
	\begin{equation}\label{1}
		\frac{\partial}{\partial t}g(x,y;t)=-2h(x,y;t),
	\end{equation}
	where $ g(x,y;t)=\frac{1}{2}\frac{\partial^2}{\partial y^i \partial y^j}F^2(x,y;t)dx^i\otimes dx^j$ is the time-dependent fundamental form of $(M,F)$, and $h(x,y;t)=h_{ij}(x,y;t)dx^i\otimes dx^j$ is a general 0-homogeneous $(0,2)$-tensor.	
	Flow \eqref{1} could also be discribed as 
	\[
	\frac{\partial}{\partial t}g_{ij}=-2h_{ij},
	\]
	in local coordinates. When $h$ is chosen to be the Akbar Zadeh's Ricci tensor, namely, $h=\mathcal{R}ic$ or $h_{ij}(x,y;t)=Ric_{ij}$, the flow \eqref{1} turns to be the Finsler Ricci flow in \cite{Lakzian:2015,Zeng-He:2019,Cheng:2022}.
	We prove the gradient estimates of positive solutions to the heat flow \eqref{2} under the Finsler geometric flow \eqref{1}. 
	\begin{theorem}\label{thm1} 
		Let $(M,F(t),\mu)$ be an $n$-dimensional compact $CD(-K,N)$ solution to the Finsler geometric flow \eqref{1}, for some $N > n$, and assume that $\vert\nabla \tau\vert\le K'$. Suppose that $-L_1g\le h\le L_1g$, $\vert \nabla h \vert \le L_2$ and $\vert \dot\nabla h \vert\leq L_3$ for some positive constants $L_1,L_2,L_3$, where  $\nabla$ and $\dot\nabla$ denote the horizontal and vertical Chern derivatives, respectively.
		Consider a positive solution $u=u(x,t)$ to the heat flow \eqref{2}. Then for any $\alpha>1$, it satisfies on $M\times [0,T]$ that 
		\begin{align}
			\frac{F^2(\nabla u)}{u^2}-\alpha\frac{u_t}{u} \le \frac{N\alpha^2}{t}&+\frac{N\alpha^2}{2}\Bigg(\frac{K-\epsilon}{\alpha-1}+\frac{K'}{2(\alpha-1)(N-n)} \notag\\
			&+\Big(1+\sqrt{2(N-n+4)}\Big)L_1+\sqrt{\frac{2}{\epsilon N}}L_2+\sqrt{\frac{8}{N}}L_3\Bigg),
		\end{align}
		for any $\epsilon>0$.
	\end{theorem}
	The conditions on $\vert\nabla \tau\vert$, $\vert \nabla h \vert$, $\vert \dot\nabla h \vert$ means that they are all bounded by some constants, which is naturally satisfied since both $\tau$ and $h$ are 0-homogeneous as are their Chern covariant derivatives. In fact, $\vert\nabla \tau\vert$, $\vert \nabla h \vert$, $\vert \dot\nabla h \vert$ are all defined on the tangent sphere bundle of the compact manifold $M$.

	The Harnack inequality is a direct corollary of the gradient estimate. That is 
	\begin{theorem}\label{thm2}
		Under the same assumptions as in Theorem \ref{thm1}. Consider two points $(x_1,t_1),\,(x_2,t_2)$ in $M^n\times(0,T)$ with $t_1<t_2$. Let $\eta(s)$ be a smooth curve connecting from $x_2$ to $x_1$, such that $\eta(1)=x_1$ and $\eta(0)=x_2$. Denote the length of $\dot{\eta}(s)$ at time $\xi(s)=(1-s)t_2+t_1$ by $F(\dot{\eta}(s))|_{\xi}$. Then, the following inequality hold.
		\begin{equation}
			u(x_1,t_1)\le u(x_2,t_2)\left(\frac{t_2}{t_1}\right)^{N{\alpha}} \exp\left\{ \int_0^1\frac{\alpha}{4}\frac{F^2(\dot{\eta}(s))|_{\xi}}{t_2-t_1}ds +\frac{N\alpha}{2}\mathcal{Q}(t_2-t_1)\right\},
		\end{equation}
		where $\alpha>1$ and $\mathcal{Q}=\mathcal{Q}(N,n,L_1,L_2,L_3,\alpha,K, K',\epsilon)$ given by 
		\begin{equation}\label{eq-Q}
			\mathcal{Q}=\frac{K-\epsilon}{\alpha-1}+\frac{K'}{2(\alpha-1)(N-n)}+\Big(1+\sqrt{2(N-n+4)}\Big)L_1+\sqrt{\frac{2}{\epsilon N}}L_2+\sqrt{\frac{8}{N}}L_3.
		\end{equation} 
	\end{theorem}
	
    An alternative condition worth considering is that on a compact manifold, the boundedness of the first derivatives $\nabla h$ and $\dot{\nabla} h$ ensures the boundedness of the evolution term $h$ itself. Consequently, in the theorems above, the condition $-L_1 g \leq h \leq L_1 g$ may be omitted. However, in this case, the resulting constants will depend on the diameter of the tangent sphere bundle $SM$.\\
    
    This manuscript is organized as follows. In Sect. 2, we give some necessary basic knowledges and lemmata in Finsler geometry, as well as the concepts and notaions in Finsler geometric flow. In Sect. 3, we give the important lemmata about the evolution of the nonlinear first and second order differential operators. In Sect. 4, we give complete proofs of the gradient estimate and the Harnack inequality.

	\section{Preliminaries}\label{sec2}
	We introduce the background in two parts. First, we recall the concept of a Finsler metric measure space. For further details, we refer the reader to \cite{Y.Shen:2016}. Subsequently, we present the related concepts concerning the Finsler geometric flow.
	\subsection{Finsler metric measure space}	
	A Finsler metric measure space $(M,F,\mu)$ is a triple consisting of a differential manifold $M$, a Finsler metric $F$ and a Borel measure $\mu$. Here, $F$ is a positive 1-homogeneous norm $F: TM \to [ 0 , \infty )$, defined on the tangent bundle, which induces the \textit{fundamental form} $g=g_{ij}dx^i\otimes dx^j$ as $g_{i j}(x, y):=\frac{1}{2} \frac{\partial^2 F^2}{\partial y^i \partial y^j}(x, y)$. Taking the derivative of $F^2$ along the fibre again yields the \textit{Cartan tensor}, Namely
	\begin{equation*}
		C_y(X,Y,Z):=C_{ijk}(y)X^iY^jZ^k:=\frac{1}{4}\frac{\partial^3 F^2(x,y)}{\partial y^i\partial y^j\partial y^k}X^iY^jZ^k,   
	\end{equation*}
	for any local vector fields $X$, $Y$, $Z$. There exists a unique almost $g$-compatible and torsion-free connection on the pull back tangent bundle $\pi^{*} TM$, known as the \textit{Chern connection}, which is defined by
	\begin{align*}
		\nabla_XY-\nabla_YX&=[X,Y]; \notag  \\
		Z(g_y(X,Y))-g_y(\nabla_ZX,Y)-g_y(X,\nabla_ZY)&=2C_y(\nabla_Zy,X,Y).
	\end{align*}
	The Chern connection
	coefficients is locally denoted by $\Gamma_{jk}^i(x, y)$ in the natural coordinate system. These coefficients
	locally induce the \textit{spray coefficients} as $G^i = \frac{1}{2}\Gamma_
	{jk}^iy^jy^k$. The spray is given by
	\begin{equation*}
		G=y^i\frac{\delta}{\delta x^i}=y^i\frac{\partial}{\partial x^i}-2G^i\frac{\partial}{\partial y^i}, 
	\end{equation*} 
	where $\frac{\delta}{\delta x^{i}}=\frac{\partial}{\partial x^{i}}-N_{i}^{j}\frac{\partial}{\partial y^{j}}$, and the nonlinear connection coefficients are locally induced
	from the spray coefficients by $N_{j}^{i}=\frac{\partial G^{i}}{\partial y^{j}}$. Traditionally, the horizontal Chern derivative is denoted by `` $\mid$ " and the vertical Chern derivative by `` ; ". For example,
	\begin{equation*}
	    w_{i|j}=\frac{\delta}{\delta x^{j}}w_{i}-\Gamma_{ij}^{k}w_{k},\quad w_{i;j}=F\frac{\partial}{\partial y^{j}}w_{i}, 
	\end{equation*}
	for any 1-form $w = w_idx^i$ on the pull-back bundle. In this manuscript, the horizontal Chern derivative is also denoted by “$\nabla$” and the vertical Chern derivative by “$\dot\nabla$”. 
	
	The \textit{Chern Riemannian curvature} $R$ is locally defined by
	\begin{equation*}
		R_{j\;kl}^{i}=\frac{\delta\Gamma_{jl}^{i}}{\delta x^{k}}-\frac{\delta\Gamma_{jk}^{i}}{\delta x^{l}}+\Gamma_{km}^{i}\Gamma_{jl}^{m}-\Gamma_{lm}^{i}\Gamma_{jk}^{m}.
	\end{equation*}
	The \textit{flag curvature} with pole $y$ is defined from $R$ as
	\begin{equation*}
		K(P,y)=K(y,u):=\frac{R_{y}(y,u,u,y)}{F^{2}(y)h_{y}(u,u)}=\frac{-R_{ijkl}(y)y^iu^jy^ku^l}{(g_{ik}(y)g_{jl}(y)-g_{il}(y)g_{jk}(y))y^iu^jy^ku^l},  
	\end{equation*}
	for any two linearly independent vectors $y$, $u$ $\in$ $T_xM \setminus \{0\}$, which span a tangent plane $\Pi_{y}=\operatorname{span}\{y,u\}$. Then the \textit{Finslerian Ricci curvature} could be given as 
	\begin{equation}
		Ric(y):=F^2(y)\sum_{i=1}^{n-1}K(y,e_i),   \notag
	\end{equation}
	where $e_1,\cdots, e_{n-1}, \frac{y}{F(y)}$ form an orthonormal basis of $T_xM$ with respect to $g_y$. 
	
	The forward distance from $p$ to $q$ is defined by
	\begin{equation}
		d_p(q):=d(p, q):=\inf _\gamma \int_0^1 F(\gamma(t), \dot{\gamma}(t)) dt,   \notag
	\end{equation}
	where the infimum is taken over all the $C^1$
	curves $\gamma$ : [0, 1] $\rightarrow$ $M$ such that $\gamma(0) = p$
	and $\gamma(1) = q$. A $C^2$ curve $\gamma$ is called a \emph{geodesic} if it locally satisfies the geodesic eqaution
	$$\ddot\gamma^i(t)+2G^i(\gamma(t),\dot \gamma(t))=0.$$
	
	The \textit{Legendre transformation} is an isomorphism $\mathcal{L}^*:T_x^*M\to T_xM$ mapping $v^*\in T_x^*M$ to a unique element $v \in T_xM$ such that $v^*( v) = F( v) $. For a differentiable function $f:M \to \mathbb{R}$, the \textit{gradient~vector} of $f$ at $x$ is defined by $\nabla( x) : = \mathcal{L} ^* ( df( x) ) \in T_xM.$ If $df( x) = 0, $ then we have $\nabla f(x) = 0.$ If $Df( x) \neq 0, $ and we can write on $M_f:=\{x\in M:df(x)\neq0\}$ that
	\begin{equation}
		\nabla f=\sum_{i,j=1}^{n}g^{ij}(x,\nabla f)\frac{\partial f}{\partial x^{i}}\frac{\partial}{\partial x^{j}}.   \notag
	\end{equation}
	In local coordinates 
	$\{x^{i}\}_{i=1}^{n}$, expressing the volume measure by $d\mu=e^{\Phi}dx^{1}\wedge\cdots\wedge dx^{n}$, then the \textit{divergence} of a smooth vector field $V$ can be written locally as
	\begin{equation}
		\mathrm{div}_{\mu}V=\sum_{i=1}^{n}(\frac{\partial V^{i}}{\partial x^{i}}+V^{i}\frac{\partial\Phi}{\partial x^{i}}).  
	\end{equation}
	The \textit{Finsler Laplacian} of a function $f$ on $M$ could also be given by
	\begin{equation}
		\Delta f:=\mathrm{div}_{\mu}(\nabla f),   \notag
	\end{equation}
	in the distributional sense. That is, 
	\begin{equation*}
		\int_{M}\phi\Delta fd\mu=-\int_{M}d\phi(\nabla f)d\mu,  
	\end{equation*}
	for any $f \in W^{1,p}(M)$ and a test function $\phi \in C_0^{\infty}(M)$.
	
	The distortion of $(M, F, \mu)$ and the \textit{S-curvature} are defined by
	\begin{equation}
		\tau(x,y):=\frac12\log\det g_{ij}(x,y)-\Phi(x), \mbox{ and } S(x,y):=\frac{d}{dt}(\tau \circ \gamma)|_{t=0},
	\end{equation}
	respectively, where $\gamma = \gamma(t)$ is a forward geodesic from $x$ with the initial tangent vector
	$\dot{\gamma}(0) = y$.

	The following concept, known as the \textit{weighted Ricci curvature}, has been proven to be equivalent to the $CD(-K,N)$ condition \cite{Ohta:2009}.  	
	\begin{definition}
		Given a unit vector $V \in T_xM$ on a Finsler metric measure space $(M,F,\mu)$ and a positive number $k$,
		the weighted Ricci curvature is defined by		
		\begin{equation}
			Ric^k(V):=\begin{cases}Ric(x,V)+\dot{S}(x,V),&\text{if }S(x,V)=0\text{ and }k=n\text{ or if }k=\infty;\\-\infty,&\text{if }S(x,V)\neq0\text{ and if }k=n;\\Ric(x,V)+\dot{S}(x,V)-\frac{S^2(x,V)}{k-n},&\text{if }n<k<\infty,\end{cases}   \notag
		\end{equation}
		where $\dot{S}(x,V)$ is the derivative along the geodesic from $x$ in the direction of $V$ .
	\end{definition}
	Based on it, we define that
	\begin{definition}
		A Finsler metric measure space satisfying $\mathrm{Ric}^{N}\geq -K$ for some $N>n$ and $K>0$ is referred to as a \textit{Finsler $CD(-K,N)$ space}.
	\end{definition}
	The term ``$CD(-K,N)$ condition" is equivalent to the weighted curvature $\mathrm{Ric}^{N}$ bounded below by a constant $-K$ in the Finsler case, as demonstrated by Ohta in \cite{Ohta:2009}. Therefore, for the sake of conciseness, we often refer to the (Finsler) metric measure space with a weighted Ricci curvature bounded from below as a (Finsler) $CD(-K,N)$ space.
	
	The following Bochner-Weitzenb\"ock formula is adopted to generalize the Li-Yau gradient estimates on Finsler $CD(-K,N)$ spaces. 
	\begin{lemma}[\cite{Ohta:2014}]\label{Bochner}
		Let $(M, F, \mu)$ be an n-dimensional Finsler metric measure space. Given $u$ $\in$ $H_{loc}^2(M) \cap C^1(M)$  with $\Delta u \in H_{loc}^1(M)$, we have on $M_u$ that 
		\begin{equation} \label{3}
			\Delta^{\nabla u}(F^2(\nabla u))=2du(\nabla^{\nabla u} (\Delta u))+2\Vert \nabla^2u\Vert^2_{HS(\nabla u)}+2 Ric^{\infty}(\nabla u).
		\end{equation}
	\end{lemma}
		
	\subsection{Finsler geometric flow}
	When considering the Finsler geometric flow and defining $h(y)=h_{ij}(y)y^iy^j$, \eqref{1} can be rewritten as 
	\begin{eqnarray}
		\frac{\partial}{\partial t}\log F=-H,
	\end{eqnarray}
	where $H=h(y)/F^2$ is a function defined on the sphere bundle $SM$. Thus, the solution to the Finsler geometric flow is a time-dependent Finsler metric $F(x,y;t)$, such that the entire metric measure space is described by $(M, F(t), \mu)$. In local coordinates, it can be directly verified that 
    \begin{eqnarray}
    	\frac{\partial}{\partial t}g^{ij}=2h^{ij},
    \end{eqnarray} 
    where $(g^{ij}(x,y;t))$ represents the inverse matrix of the fundamental form, and $h^{ij}(x,y;t)=g^{ik}h_{kl}g^{lj}$. 
    
    Given a fixed local vector field $V$ on $M$, the expression $h(x,V;t)=h_{ij}(V)dx^i\otimes dx^j$ defines a time-dependent symmetric bilinear form, that is,
    \begin{eqnarray}\label{a2.6}
    	h_V(W_1,W_2)=h_{ij}(V)W_1^iW_2^j,
    \end{eqnarray}
    for any two vector fields $W_1$ and $W_2$. Obviously, $h_V(V,V)=h(V)$. 
    When $h_V(\cdot,\cdot)$ acts on a (0,2)-type tensor $\Omega$, we simply denote it by $h_V(\Omega)$.
    
    The form $h_V(\cdot,\cdot)$ can be pull back to a symmetrical bilinear $(2,0)$-type tensor, denoted by $h^{\flat}_V$. Namely, 
    \begin{eqnarray}
    	h^{\flat}_V(\alpha,\beta)=h^{ij}(V)\alpha_i\beta_j,
    \end{eqnarray}
    for any two forms $\alpha$ and $\beta$ on $M$, where $h^{ij}(V)=g^{ik}(V)h_{kl}(V)g^{lj}(V)$. We denote the vector field $h^{\flat}_V(\alpha,\cdot)$ by $\alpha(h^{\flat}_V)$ for short.
    
    The horizontal Chern derivative of $h_{V}$ with respect to the reference vector $V$ is expressed as 
    \begin{eqnarray}
    	(\nabla^{V}h)_{V}(W_1,W_2,W_3)=h_{ij|k}(V)W_1^iW_2^jW_3^k,
    \end{eqnarray}
    for any vector fields $W_1,W_2,W_3$, where $\nabla^V$ indicates that the horizontal Chern derivative is taken with respect to the reference vector field $V$, and $(\cdot)_V$ signifies that the tensor is evaluated at $y=\nabla f$. The trace of $(\nabla^{V}h)_{V}$ with respect to $V$ is a 1-form denoted by 
    \begin{eqnarray}\label{trnablah}
    	(\mathrm{tr}_{V}\nabla^Vh)_{V}=g^{ik}(V)h_{ij|k}(V)dx^j=h^i_{\,j|i}(V)dx^j.
    \end{eqnarray}
     It is important to note the indices used for the trace in \eqref{trnablah}. It should to be the indices $i,k$ in $h_{ij|k}$, as the notation should reflect $\nabla^V(\mathrm{tr}_Vh)_V$ when considering $i,j$ in $h_{ij|k}$.
    
    Similarly, the vertical Chern derivative of $h^{\flat}_V$ with respect to $V$ is denoted by $(\dot\nabla h^{\flat})_V$ and is defined as  
    \begin{eqnarray}\label{a2.10}
    	(\dot\nabla h^{\flat})_V(\alpha,\beta,\cdot)=h^{ij}_{\,\,;k}(V)\alpha_i\beta_j\delta y^k(V),
    \end{eqnarray} 
    for any 1-forms $\alpha$, $\beta$, where $\delta y^k(V)=dy^k+N^k_l(V)dx^l$. 
    The (1,1)-type tensor $(\dot\nabla h^{\flat})_V(\alpha,\cdot,\cdot)$ can be succinctly denoted by $\alpha((\dot\nabla h^{\flat})_V)$.                                                     
	
	\section{Derivative operator exchange formulae and some lemmata}\label{sec3}
	In this section, we present several lemmata concerning functions on a time-dependent Finsler metric measure space, as well as solutions to the heat flow \eqref{2} under the Finsler geometric flow \eqref{1}.  
	
	\begin{lemma}\label{lem1}
		Let $(M,F(t),\mu)$ be a solution to the Finsler geometric flow \eqref{1}. For any $f\in C^1(M)\cap C^1([0,T])$, it satisfies that
		\begin{equation}\label{eq-lem1}
			\partial_t(F^2(\nabla f))=2h(\nabla f)+2df_{t}(\nabla f),
		\end{equation}
		where $h(\nabla f)=h_{\nabla f}(\nabla f, \nabla f)=h_{ij}(\nabla f)f^if^j$.
	\end{lemma}
	\begin{proof}
		 The lemma follows from a direct computation. That is,
		\begin{align*}
			\partial_t(F^2(\nabla f))=&\partial_t(g^{ij}(\nabla f)f_if_j)\notag\\
			=&(\partial_tg^{ij})(\nabla f)f_if_j+\frac{\partial g^{ij}}{\partial y^k}\frac{\partial f^k}{\partial t}f_if_j+2g^{ij}(\nabla f)\frac{\partial f_i}{\partial t}f_j\notag\\
			=& 2h^{ij}(\nabla f)f_if_j-2C_{\;k}^{ij}(\nabla f)\frac{\partial f^k}{\partial t}f_if_j+2g^{ij}(\nabla f)\frac{\partial f_i}{\partial t}f_j\notag\\
			=& 2h^{ij}(\nabla f)(\nabla f, \nabla f)+2df_{t}(\nabla f).
		\end{align*}
	\end{proof}

	Since the solutions to \eqref{2} on $(M,F(t),\mu)$ lack sufficient regularity, the following derivative operator exchange formulae are interpreted in a weak sense.
	\begin{lemma}\label{lem2}
		Let $(M,F(t),\mu)$ be a solution to \eqref{1}. For any $f\in H^1([0,T],H^1(M))\cap L^2([0,T],H^1(M))$,
		the following two operators exchange formulae are satisfied on $M_f$, hence, are satisfied in the distributional sense. 
		\begin{itemize}
			\item[i)] $[\nabla^{\nabla f},\partial_t]f=-2df(h^{\flat}_{\nabla f})$, where $df(h^{\flat}_{\nabla f})=h^{ij}(\nabla f)f_i\frac{\partial}{\partial x^i}$ is a vector field on $M$.
			\item[ii)] $[\Delta^{\nabla f},\partial _t]f=-2\mathcal{J}(f,h,F,\mu)$, where 
			\begin{eqnarray}\label{J-g}
				\begin{split}
					\mathcal{J}(f,h,F,\mu)=&h_{\nabla f}(\nabla^2f)+(\mathrm{tr}_{\nabla f}\nabla^{\nabla f}h)_{\nabla f}(\nabla f)\\
					&+\frac1Fdf((\dot\nabla h^{\flat})_{\nabla f})(d\nabla f)-h_{\nabla f}(\nabla f,\nabla^{\nabla f}\tau).
				\end{split}
			\end{eqnarray}
		\end{itemize}    
	\end{lemma}
	\begin{proof}
		The first formula can be direct verified as 
		\begin{align*}
			\nabla^{\nabla f}(f_t)=&g^{ij}(\nabla f)\frac{\partial f_t}{\partial x^i} \frac{\partial}{\partial x^j}\notag\\
			=& g^{ij}(\nabla f)\partial_t(f_i) \frac{\partial}{\partial x^j}\notag\\
			=& \partial_t\left(g^{ij}(\nabla f)f_i\frac{\partial}{\partial x^j}\right)-\partial_t\left(g^{ij}(\nabla f)\right)f_i\frac{\partial}{\partial x^j}\notag\\
			=& \partial_t(\nabla f)-2h^{ij}(\nabla f)f_i\frac{\partial}{\partial x^j}.
		\end{align*}
		
		To prove the second one, we choose a test function $\phi \in H_0^1([0,T]\times M)$, so that
		\begin{align*}
			\int_{0}^T\int_M \phi \Delta^{\nabla f}f_t d\mu dt=&-\int_0^T\int _M d\phi (\nabla^{\nabla f}f_t)d\mu dt\notag\\
			=& -\int_0^T\int_M d\phi\left[\partial_t(\nabla f)-2df(h^{\flat}_{\nabla f})\right]d\mu dt\notag\\
			=&\int _0^T\int_M d(\phi_t)(\nabla f)dmdt+2\int_0^T\int_Mh_{\nabla f}(\nabla^{\nabla f}\phi,\nabla f)d\mu dt.
		\end{align*}
		Then the formula \textit{ii)} is followed from the next lemma.
	\end{proof}
	
	Before completing the proof, we provide an explanation of the geometric meaning of $\mathcal{J}$. In local coordinates, it can be rewritten as 
	\begin{eqnarray}\label{J-l}
		\mathcal{J}(f,h,F,\mu)=h^{ij}(\nabla f)f_{j|i}+h_{\; |i}^{ij}(\nabla f)f_j+\frac{1}{F}h_{\; ;k}^{ij}f_{\;|i}^kf_j-h^{ij}(\nabla f)f_i\tau_{|j},
	\end{eqnarray} 
	where $\tau$ represents the distortion of $(M,F,\mu)$. The correspondence between \eqref{J-g} and \eqref{J-l} is established in \eqref{a2.6}-\eqref{a2.10}. The quantity $\mathcal{J}$ essentially captures the difference in the geometric flow evolution of the second-order quasilinear elliptic operator $\Delta^{\nabla f}$ at each time slice, as well as the evolution of $\Delta^{\nabla f}$ along the flow. It not only depends on the geometry of the horizontal and vertical bundles but is also closely related to the Borel measure $\mu$, since the nonlinear Laplacian is highly dependent on it.
	\begin{lemma}\label{lem3}
		Let $(M,F(t),\mu)$ be a solution to \eqref{1}, and let $f\in C^1([0,T])\cap C^2(M)$ be a function on $M\times[0,T]$.  Then it satisfies that 
		\begin{equation*}
			h_{\nabla f}(\nabla f,\nabla^{\nabla f}\phi )=-\phi \mathcal{J}+\mathrm{div}(\phi df(h^{\flat}_{\nabla f})),
		\end{equation*}
		for any $\phi \in H_0^1(M\times [0,T])$.
	\end{lemma}
	\begin{proof}
		According to the definition of $h_{\nabla f}(\cdot,\cdot)$, we have that
		\begin{align}\label{a3.2}
			\int_{0}^T\int_M &h_{\nabla f}(\nabla f,\nabla^{\nabla f}\phi)d\mu dt=\int_{0}^T\int_M\phi_ih^{ij}(\nabla f)f_j e^{\Phi}dxdt\notag\\
			=& \int_{0}^T\int_M \frac{\partial }{\partial x^i}\left(\phi h^{ij}(\nabla f)f_j e^{\Phi}\right)e^{-\Phi}d\mu dt- \int_{0}^T\int_M\phi \frac{\partial}{\partial x^i}\left(h^{ij}(\nabla f)f_je^{\Phi}\right)dxdt
			\notag\\
			=& \int_{0}^T\int_M c(\phi df(h^{\flat}_{\nabla f}))d\mu dt-\mathcal{A},
		\end{align}
		where we denote the last term on the right-hand side (RHS) of \eqref{a3.2} by $\mathcal{A}$. That is, 
		\begin{align*}
			\mathcal{A}=& \int_{0}^T\int_M\phi\left[ \left(\frac{\partial h^{ij}}{\partial x^i}(\nabla f)+\frac{\partial h^{ij}}{\partial y^k}\frac{\partial f^k}{\partial x^i}\right)f_j+h^{ij}(\nabla f)\frac{\partial f_j}{\partial x^i}+h^{ij}(\nabla f)f_j\Phi_i\right]d\mu dt\notag\\
			=&\int_{0}^T\int_M\phi \Bigg\{  \left[\frac{\delta h^{ij}}{\delta x^i}(\nabla f)+N_i^k(\nabla f)\frac{\partial h^{ij}}{\partial y^k}(\nabla f)+ \frac{\partial h^{ij}}{\partial y^k}\left(\frac{\delta f^k}{\delta x^i}+N_{i}^l(\nabla f)\frac{\partial g^{km}}{\partial y^l}f_m\right)\right]f_j\notag \\
			&\quad\quad\quad\quad\quad+h^{ij}(\nabla f)\frac{\partial f_j}{\partial x^i}+h^{ij}(\nabla f)f_j\Phi_i\Bigg\}d\mu dt\notag\\
			=&\int_{0}^T\int_M\phi\Bigg\{ \left[h_{\;|i}^{ij}(\nabla f)-\Gamma_{il}^ih^{lj}-\Gamma_{il}^jh^{il}+\frac{\partial h^{ij}}{\partial y^k}\left(\frac{\delta f^k}{\delta x^i}+N_i^l(\nabla f)\right)\right]f_j\notag\\
			&\quad\quad\quad\quad\quad+h^{ij}(\nabla f)\frac{f_j}{x^i}+h^{ij}(\nabla f)f_j\Phi_i\Bigg\} d\mu dt\notag\\
			=&\int_{0}^T\int_M\phi\left\{ h_{\;|i}^{ij}(\nabla f)f_j+h^{ij}(\nabla f)f_{j|i}-h^{ij}(\nabla f)f_j\tau_{|i}+\frac{1}{F}h_{\; ;k}^{ij}f_{\; |i}^kf_j\right\}d\mu dt\notag\\
			=& \int_0^T\int _M \phi \mathcal{J}d\mu dt.
		\end{align*}
	\end{proof}
	
	The distributional sense of \textit{ii)} in Lemma \ref{lem2} is equal to that
	\begin{eqnarray}\label{lem2ii-2}
		\int_0^T\int_M\left[d\phi_t(\nabla f)+d\phi(\nabla^{\nabla f}f_t)\right]d\mu dt-2\int_0^T\int_M\phi\mathcal{J} d\mu dt=0,
	\end{eqnarray}
	for any $\phi\in C^{\infty}_c(M\times [0,T])$.
	
	\begin{remark}\label{rmk1}
		The omission in \cite{Lakzian:2015,Zeng-He:2019,Zeng:2020,Zeng:2021,Cheng:2022,Azami:2023} is that the effect of the measure on the operator exchange given above is not cinsidered. This oversight manifests appeared in the derivation of several key equations, including: equation (48) in the proof of Lemma 4.1 in \cite{Lakzian:2015}; equation (3.11) in the proof of Lemma 3.3 in \cite{Zeng-He:2019}; equation (3.15) in the proof of Lemma 3.4 in \cite{Zeng:2020}; equation (4.7) in the proof of Lemma 4.3 in \cite{Zeng:2021}; the calculation of $\mathcal{A}$ in the proof of Lemma 3.2 in \cite{Cheng:2022} and the calculation of $\mathcal{A}$ in the proof of Lemma 3.4 in \cite{Azami:2023}, among others. 
	\end{remark}
	\begin{remark}\label{rmk2}	
	Additionally, the requirement of the vanishing of the S-curvature in \cite{Lakzian:2015,Zeng-He:2019,Zeng:2020,Zeng:2021} is  overly restrictive, as the S-curvature of a solution $(M,F(t),\mu)$ to the Finsler geometric flow (or the Finsler Ricci flow) inherently depends on the time parameter $t$. In fact, the results in \cite{Zeng-He:2019,Zeng:2020,Zeng:2021,Cheng:2022,Azami:2023} can only be established under the additional assumption that the distortion $\tau$ is Chern horizontal invariant.\\ 
	\end{remark}
	
	
	We now turn our attention to the solution $u=u(x,t)$ of the heat equation \eqref{2} under the Finsler geometric flow. In this manuscript, we focus exclusively on global solutions to the Finsler heat equation \eqref{2}. More general nonlinear parabolic equations will be investigated in subsequent papers.

	Let $u$ be a positive solution to Finsler heat equation $\partial_tu=\Delta u$ on $(M,F(t),\mu)$, that is $u=u(s,t)\in L^2([0,T],H_0^1(M))\bigcap H^1([0,T],H^{-1}(M))$ satisfies
	\begin{equation}
		\int _M \phi \partial_t ud\mu=-\int_Md\phi(\nabla u)d\mu,
	\end{equation}
	 for any $\phi \in C_c^{\infty}(M)$ and for almost all $t\in [0,T]$. 	
	Let $f(x,t)=\log u(x,t)$, one can easily find that 
	\begin{lemma}\label{lem4}
		The function $f$ satisfies that
		\begin{eqnarray}\label{eq-lem4}
			\int_0^T\int_M(d\phi(\nabla f)-\phi_tf)d\mu dt=\int_0^T\int_M\phi F^2(\nabla f)d\mu dt,
		\end{eqnarray}
		for any $\phi\in C^{\infty}_c(M\times[0,T])$.
	\end{lemma}
	Define 
	\begin{eqnarray}
		\sigma(x,t):=tf_t(x,t),\quad\mbox{and}\quad\mathcal{F}(x,t)=tF^2(\nabla f)(x,t)-\alpha\sigma(x,t),
	\end{eqnarray} 
	where $\alpha>1$ is a constant. A straightforward calculation yields the following Lemma.	
	\begin{lemma}\label{lem5}
		The function $\sigma$ satisfies the parabolic equation in the sense of distributions as
		\begin{equation}\label{eq-3.5}
			\int_0^T\int_M\left[\phi_t\sigma-d\phi(\nabla^{\nabla f}\sigma)+\frac{\phi\sigma}{t}+2\phi d\sigma(\nabla f)\right]d\mu dt=-2\int_0^T\int_Mt\phi[h(\nabla f)+\mathcal{J}]d\mu dt,
		\end{equation}
		for any $\phi \in C_c^{\infty}(M)$ and for almost all $t\in [0,T]$.
	\end{lemma}
	\begin{proof}
		According to the definition of $\sigma$, the left-hand side (LHS) of \eqref{eq-3.5} is equal to 
		\begin{align*}
			LHS\,\, of\,\, \eqref{eq-3.5}&=\int_0^T\int_M[t\phi_tf_t-td\phi(\nabla^{\nabla f}f_t)+\phi f_t+2t\phi df(\nabla^{\nabla f}f_t)]\notag\\
			&=\int_0^T\int_M[(t\phi)_tf_t-d(t\phi)(\nabla^{\nabla f}f_t)+2(t\phi)df(\nabla^{\nabla f}f_t)]d\mu dt,
		\end{align*}
		for any $\phi \in C_c^{\infty}(M)$ and for almost all $t\in [0,T]$. Taking considering of Lemma \ref{lem4} by replacing the test function $\phi$ by $(t\phi)_t$, we have
		\begin{align}
			LHS\,\, of\,\, \eqref{eq-3.5}&=\int_0^T\int_M[(t\phi)_tF^2(\nabla f)-d(t\phi)_t(\nabla f)-d(t\phi)(\nabla^{\nabla f}f_t)+2(t\phi)df(\nabla^{\nabla f}f_t)]d\mu dt \notag\\
			&=\int_0^T\int_M[(t\phi)_tF^2(\nabla f)-2t\phi\mathcal{J}+2(t\phi)df(\nabla^{\nabla f}f_t)]d\mu dt \notag\\
			&=-2\int_0^T\int_Mt\phi[h(\nabla f)+\mathcal{J}]d\mu dt,
		\end{align}
		where the second equality follows from \eqref{lem2ii-2}, and the last equality holds according to the distributional sense of \eqref{eq-lem1} with the test function $t\phi$.		
	\end{proof}	
	From this, we can deduce the equation satisfied by $\mathcal{F}$. Specifically, 
	\begin{lemma}\label{lem6}
		$\mathcal{F}(x,t)$ satisfies that
		\begin{align}\label{eq-3.6}
			\int_0^T&\int_M\left[\phi_t\mathcal{F}-d\phi(\nabla^{\nabla f}\mathcal{F})+\phi\left(2d\mathcal{F}(\nabla f)+\frac{\mathcal{F}}{t}\right)\right]d\mu dt\notag\\ 
			&=2\int_0^T\int_Mt\phi\left[ (\alpha-1)h(\nabla f)+Ric^{\infty}(\nabla f)+\Vert \nabla^2f \Vert ^2_{HS(\nabla f)}+\alpha\mathcal{J}\right]d\mu dt,
		\end{align}
		for any $\phi\in C_c^{\infty}(M\times [0,T])$.
	\end{lemma}
	\begin{proof}
		As demonstrated in the proof of Lemma \ref{lem5}, we still present the calculations in weak sense.
		\begin{align}\label{mBC}
			LHS\,\,of \,\,\eqref{eq-3.6}=&\int_0^T\int_M\left[(t\phi)_tF^2(\nabla f)-d(t\phi)(\nabla^{\nabla f}F^2(\nabla f))+2(t\phi)df(\nabla^{\nabla f}F^2(\nabla f))\right]d\mu dt \notag\\
			&+\alpha\int_0^T\int_M\left[d\phi(\nabla^{\nabla f}\sigma)-\phi_t\sigma-2\phi d\sigma(\nabla f)-\phi\frac{\sigma}{t} \right]d\mu dt \notag\\
			=:&\mathcal{B}+\mathcal{C},
		\end{align}
		where 
		\begin{align}\label{mC}
			\mathcal{C}=&\alpha\int_0^T\int_M\left[d\phi(\nabla^{\nabla f}\sigma)-\phi_t\sigma-2\phi d\sigma(\nabla f)-\phi\frac{\sigma}{t} \right]d\mu dt \notag\\
			=&2\alpha\int_0^T\int_Mt\phi[h(\nabla f)+\mathcal{J}]d\mu dt,
		\end{align} 
		by Lemma \ref{lem5}, and 
		\begin{align}
			\mathcal{B}=&\int_0^T\int_M\left[(t\phi)_tF^2(\nabla f)-d(t\phi)(\nabla^{\nabla f}F^2(\nabla f))+2(t\phi)df(\nabla^{\nabla f}F^2(\nabla f))\right]d\mu dt \notag\\
			=&\int_0^T\int_M\left[(t\phi)_tF^2(\nabla f)+2(t\phi)df(\nabla^{\nabla f}F^2(\nabla f))\right]d\mu dt\notag\\
			&-\int_0^T\int_{M_f}d(t\phi)(\nabla^{\nabla f}F^2(\nabla f))d\mu dt,
		\end{align}
		according to the Sard Theorem. Recalling the Bochner-Weitzenb\"ock formula in Lemma \ref{Bochner}, we can derive that
		\begin{align}\label{mB}
			\mathcal{B}=&\int_0^T\int_M\left[(t\phi)_tF^2(\nabla f)+2(t\phi)df(\nabla^{\nabla f}F^2(\nabla f))\right]d\mu dt\notag\\
			&+2\int_0^T\int_{M_f}(t\phi)[df(\nabla^{\nabla f}(\Delta f))+\|\nabla^2f\|_{HS(\nabla f)}^2+Ric^{\infty}(\nabla f)]d\mu dt\notag\\
			=&\int_0^T\int_M(t\phi)_tF^2(\nabla f)d\mu dt+2\int_0^T\int_{M_f}(t\phi)[df_t(\nabla f)+\|\nabla^2f\|_{HS(\nabla f)}^2+Ric^{\infty}(\nabla f)]d\mu dt\notag\\
			=&-2\int_0^T\int_M(t\phi)h(\nabla f)d\mu dt+2\int_0^T\int_{M_f}(t\phi)[\|\nabla^2f\|_{HS(\nabla f)}^2+Ric^{\infty}(\nabla f)]d\mu dt.
		\end{align}
		where the second equality follows from the strong sense of \eqref{eq-lem4} on $M_f$, and the last equality holds according to the distributional sense of \eqref{eq-lem1} again.
		
		Plugging \eqref{mC} and \eqref{mB} into \eqref{mBC} finishes the proof.
	\end{proof}
	
	With these preparations, we now have all the necessary components to proceed with the proofs of the main theorems.
	
	\section{Global gradient estimates under compact Finsler $CD(-K,N)$ geometric flow}\label{sec4}
	In this section, we present the proofs of Theorems \ref{thm1} and \ref{thm2}.

	\begin{proof}[Proof of Theorem \ref{thm1}]
		One could deduce the following five lower bounds estimates by using the Cauchy-Schwartz and Young's inequalities, with constants $\epsilon_1$,$\epsilon_2$,$\epsilon_3$ and $\epsilon_4$ to be determined.
		\begin{align}
			 (\alpha-1)h(\nabla f)&=(\alpha-1)h^{ij}(\nabla f)f_if_j \ge -(\alpha-1)L_1F^2(\nabla f),\\
			 \alpha h_{\nabla f}(\nabla^2f)&=\alpha h^{ij}(\nabla f)f_{j|i} \ge -\frac{\alpha^2}{4\epsilon_1}\Vert h_{\nabla f}\Vert_{HS(\nabla f)}^2-\epsilon_1\Vert \nabla^2 f\Vert_{HS(\nabla f)}^2,\\
			 -\alpha h_V(\nabla f,\nabla^{\nabla f}\tau)&=-\frac{\alpha}{2} h^{ij}(\nabla f)(f_i\tau_{|j}+f_j\tau_{|i})\notag\\
			 &\ge -\frac{\alpha^2}{4\epsilon_2}\Vert h_{\nabla f}\Vert_{HS(\nabla f)}^2-\frac{\epsilon_2}{2}\left(F^2(\nabla f)F^2(\nabla^{\nabla f} \tau)+S^2(\nabla f)\right),\\ 
			 \alpha(\mathrm{tr}_{\nabla f}\nabla^{\nabla f}h)_{\nabla f}&(\nabla f)=\alpha h_{\;|i}^{ij}(\nabla f)f_j\ge -\frac{\alpha^2}{4\epsilon_3}F^{*2}((tr_{\nabla f}\nabla^{\nabla f}h)_{\nabla f})-\epsilon_3F^2(\nabla f),
		\end{align}
        and
		\begin{align}
			\frac{\alpha}Fdf((\dot\nabla h^{\flat})_{\nabla f})(d\nabla f)&=\frac{\alpha}{F}h_{\; ;k}^{ij}f_{\;|i}^kf_j \notag\\
			& \ge -\epsilon_4\Vert d\nabla f\Vert_{HS(\nabla f)}^2-\frac{\alpha^2}{4\epsilon_4}\left\Vert\frac{df(\dot\nabla h^{\flat}_{\nabla f})}{F}\right\Vert_{HS(\nabla f)}^2 \notag\\
			&\ge -\epsilon_4\Vert \nabla^2 f\Vert_{HS(\nabla f)}^2-\frac{\alpha^2}{4\epsilon_4}\Vert\dot\nabla h^{\flat}_{\nabla f}\Vert_{HS(\nabla f)}^2,
		\end{align}
		where 
		$\dot\nabla h^{\flat}_{\nabla f}=h^{ij}_{\,\,;k}(\nabla f)\frac{\delta}{\delta x^i}(\nabla f)\otimes\frac{\delta}{\delta x^j}(\nabla f)\otimes \delta y^k(\nabla f)$ is a $(2,1)$-tensor whose Hilbert-Schmidt norm is given by
		\begin{align}
			\Vert\dot\nabla h^{\flat}_{\nabla f}\Vert_{HS(\nabla f)}^2=g_{ij}(\nabla f) g_{kl}(\nabla f)g^{mh}(\nabla f)h_{\,\,;m}^{ik}(\nabla f)h_{\,\,;h}^{jl}(\nabla f).
		\end{align}
		Noticing that $\Vert \nabla^2 f\Vert_{HS(\nabla f)}^2=\Vert d\nabla f\Vert_{HS(\nabla f)}^2$,
		thus, according to Lemma \ref{lem6}, we have for any positive test function $\phi\in C_c^{\infty}(M\times [0,T])$ that
		\begin{align}\label{equ1}
			\int_0^T&\int_M\left[\phi_t\mathcal{F}-d\phi(\nabla^{\nabla f}\mathcal{F})+\phi\left(2d\mathcal{F}(\nabla f)+\frac{\mathcal{F}}{t}\right)\right]d\mu dt \notag\\
			\ge& 2\int_0^T\int_Mt\phi\Big\{ \left[-(\alpha-1)L_1-\epsilon_3-\frac{\epsilon_2}{2}F^2(\nabla^{\nabla f}\nabla \tau)\right]F^2(\nabla f)\notag\\
			&\quad\quad\quad\quad\quad\,\,+Ric^{\infty}(\nabla f)+(1-\epsilon_1-\epsilon_4)\Vert\nabla^2 f\Vert_{HS(\nabla f)}^2\notag\\
			&\quad\quad\quad\quad\quad\,\,-\frac{\epsilon_2}{2}S^2(\nabla f)-\left(\frac{\alpha^2}{4\epsilon_1}+\frac{\alpha^2}{4\epsilon_2}\right)\Vert h_{\nabla f}\Vert_{HS(\nabla f)}^2\notag\\
			&\quad\quad\quad\quad\quad\,\,-\frac{\alpha^2}{4\epsilon_3}F^2(tr_{\nabla f}\nabla^{\nabla f}h_{\nabla f})-\frac{\alpha^2}{4\epsilon_4}\Vert\dot\nabla h^{\flat}_{\nabla f}\Vert_{HS(\nabla f)}^2\Big\}d\mu dt.
		\end{align}
		Since it satisfies on $M_f$ that
		\begin{equation}
			\Vert \nabla^2f\Vert_{HS(\nabla f)}^2\ge \frac{1}{n}(tr_{\nabla f}\nabla^2 f)^2=\frac{1}{n}(\Delta f+S^2(\nabla f))\ge \frac{(\Delta f)^2}{N}-\frac{S^2(\nabla f)}{N-n},
		\end{equation}
		we can get from \eqref{equ1} that
		\begin{align}\label{eq-4.9}
			\int_0^T&\int_M\left[\phi_t\mathcal{F}-d\phi(\nabla^{\nabla f}\mathcal{F})+\phi\left(2d\mathcal{F}(\nabla f)+\frac{\mathcal{F}}{t}\right)\right]d\mu dt \notag\\
			\ge&2\int_0^T\int_Mt\phi\Big\{-\left[(\alpha-1)L_1+\epsilon_3+\frac{\epsilon_2}{2}F^2(\nabla^{\nabla f}\nabla \tau)\right]F^2(\nabla f)+Ric^{\infty}(\nabla f)\notag\\
			&\quad\quad\quad\quad\quad\,\,\,\,\,+(1-\epsilon_1-\epsilon_4)\frac{(\Delta f)^2}{N}
			-\left(\frac{1-\epsilon_1-\epsilon_4}{N-n}+\frac{\epsilon_2}{2}\right)S^2(\nabla f)\notag\\
			&\quad\quad\quad\quad\quad\,\,\,\,\,-\left(\frac{\alpha^2}{4\epsilon_1}+\frac{\alpha^2}{4\epsilon_2}\right)\Vert h_{\nabla f}\Vert_{HS(\nabla f)}^2-\frac{\alpha^2}{4\epsilon_3}F^2(tr_{\nabla f}\nabla^{\nabla f}h_{\nabla f})\notag\\
			&\quad\quad\quad\quad\quad\,\,\,\,\,-\frac{\alpha^2}{4\epsilon_4}\Vert \dot\nabla h^{\flat}_{\nabla f}\Vert_{HS(\nabla f)}^2\Big\}d\mu dt,
		\end{align}
		for any positive test function $\phi\in C_c^{\infty}(M\times [0,T])$.
		
		On the other hand, it follows from Lemma \ref{lem4} that	
		\begin{align}\label{equ2}
			\Delta f=-F^2(\nabla f)+f_t=-\frac{1}{\alpha}\left(\frac{\mathcal{F}}{t}+(\alpha-1)F^2(\nabla f)\right),
		\end{align}
		on $M_f$. 
		Thus, we can deduce from \eqref{eq-4.9} and the bounded assumptions of $\nabla \tau$ and $h$ that
		\begin{equation}\label{eq-4.11}
			\int_0^T\int_M\left[\phi_t\mathcal{F}-d\phi(\nabla^{\nabla f}\mathcal{F})+\phi\left(2d\mathcal{F}(\nabla f)+\frac{\mathcal{F}}{t}\right)\right]d\mu dt \ge 2\int_0^T\int_Mt\phi(A+B+C)d\mu dt,
		\end{equation}
		with 
		\begin{align}
				A&=-\left[(\alpha-1)L_1+\epsilon_3+\frac{\epsilon_2}{2}K'\right]F^2(\nabla f)+\frac{(1-\epsilon_1-\epsilon_4)}{N\alpha^2}\left(\frac{\mathcal{F}}{t}+(\alpha-1)F^2(\nabla f)\right)^2, \\
				B&= Ric^{\infty}-\left(\frac{1-\epsilon_1-\epsilon_4}{N-n}+\frac{\epsilon_2}{2}\right)S^2(\nabla f), \\
				C&= -\left(\frac{\alpha^2}{4\epsilon_1}+\frac{\alpha^2}{4\epsilon_2}\right)\Vert h_{\nabla f}\Vert_{HS(\nabla f)}^2-\frac{\alpha^2}{4\epsilon_3}F^2(tr_{\nabla f}\nabla^{\nabla f}h_{\nabla f})-\frac{\alpha^2}{4\epsilon_4}\Vert \dot\nabla h^{\flat}_{\nabla f}\Vert_{HS(\nabla f)}^2.
		\end{align}
		Setting $\epsilon_1=\epsilon_4=\frac{1}{4}$ and $\epsilon_2=\frac{1}{N-n}$ yields that 
		\begin{eqnarray}
			B=Ric^{\infty}(\nabla f)-\frac{S^2(\nabla f)}{N-n}=Ric^N(\nabla f).
		\end{eqnarray}
        In this situation, by assumptions $\nabla h$ and $\dot\nabla h$,  one may find that
        \begin{equation}\label{eq-C>}
			C\ge -\frac{n\alpha^2}{4}(4+N-n)L_1^2-\frac{\alpha^2}{4\epsilon_3}L_2^2-\alpha^2L_3^2.
		\end{equation}		
		To handle the term $A$, we set $\frac{ tF^2(\nabla f)}{\mathcal{F}(x,t)}=\mu\ge0$, otherwise $\mathcal{F}(x,t)< 0$, in which case the assertion becomes trivial.		
		Now $A$ could be written as 
		\begin{equation}\label{eq-4.17}
			A=-\left[(\alpha-1)L_1+\epsilon_3+\frac{K'}{2(N-n)}\right]\mu\frac{\mathcal{F}}{t}+\frac{(1+(\alpha-1)\mu)^2}{2N\alpha^2}\left(\frac{\mathcal{F}}{t}\right)^2.
		\end{equation}

		Since $Ric^N(\nabla f)\ge -KF^2(\nabla f)$ on a $CD(-K,N)$ space, \eqref{eq-4.11}-\eqref{eq-4.17} provide that 
		\begin{align}\label{equ3}
			&\int_0^T\int_M\left[\phi_t\mathcal{F}-d\phi(\nabla^{\nabla f}\mathcal{F})+2\phi d\mathcal{F}(\nabla f)\right]d\mu dt \notag\\
			 \ge& -\int_0^T\int_M\phi\frac{\mathcal{F}}{t}d\mu dt+2\int_0^T\int_Mt\phi\Bigg\{-\left[(\alpha-1)L_1+K+\epsilon_3+\frac{K'}{2(N-n)}\right]\mu \frac{\mathcal{F}}{t} \notag\\
			 &\quad\quad\quad\quad\quad\quad\quad\quad\quad\quad\quad\quad\quad\quad\quad\,\,\,+\frac{(1+(\alpha-1)\mu)^2}{2N\alpha^2}\left(\frac{\mathcal{F}}{t}\right)^2 +C\Bigg\}d\mu dt,
		\end{align}
		for any positive $\phi\in C_c^{\infty}(M\times [0,T])$, where $C$ has a lower bound given in \eqref{eq-C>}.

		Now, fix an arbitrary time $t\in(0,T]$, and let $\mathcal{F}$ attain its maximum at $(x_0,t_0)\in M\times [0,t]$. Since $\mathcal{F}(x_0,t_0)>0$ (othertwise, the assertion is trivial), we have $t_0\in (0,t]$ and $\mathcal{F}$ is positive on a neighborhood $U$ of $(x_0,t_0)$.
		Consequently, the RHS of \eqref{equ3} is nonpositive on $U$. If this were not the case, $\mathcal{F}$ would act as a local weak subsolution to 
		\begin{equation*}
			(\Delta^{\nabla f}+2df\cdot \nabla^{\nabla f}-\partial_t)\mathcal{F}\ge 0,
		\end{equation*}
		implying that its maximum must occur on the boundary of $U$. This would contradict the fact that $(x_0,t_0)$ is an interior point of $U$. Therefore, we get at $(x_0,t_0)$ that 
		\begin{align*}
			\mathcal{F}^2\le \frac{N\alpha^2t_0}{(1+(\alpha-1)\mu)^2}\left\{2\left[(\alpha-1)L_1+K+\epsilon_3+\frac{K'}{2(N-n)}\right]\mu+\frac{1}{t_0}\right\}\mathcal{F}-\frac{2N\alpha^2t_0^2C}{(1+(\alpha-1)\mu)^2}.
		\end{align*}
		
		Recall that $a\le a_2+\sqrt{a_1}$, if $a^2\le a_1+a_2a$ with $a,a_1,a_2>0$.
		It implies that
		\begin{align}\label{equ4}
			\mathcal{F}\le \frac{N\alpha^2t_0}{(1+(\alpha-1)\mu)^2}\left[2\left((\alpha-1)L_1+K+\epsilon_3+\frac{K'}{2(N-n)}\right)\mu+\frac{1}{t_0}\right]+\frac{\sqrt{2N}\alpha t_0}{1+(\alpha-1)\mu}\sqrt{-C}.
		\end{align}
		Since $\frac{1}{(1+(\alpha-1)\mu)^2}\le 1$ and 
		\[
		\frac{2\mu}{(1+(\alpha-1)\mu)^2} \le \frac{2\mu}{4(\alpha-1)\mu}=\frac{1}{2(\alpha-1)},
		\]
		for any nonnegative $\mu$, \eqref{equ4} implies the following by setting $\epsilon_3=\epsilon$,		
		\begin{align}
			\mathcal{F}\le &\frac{N\alpha^2t_0}{2(\alpha-1)}\left((\alpha-1)L_1+K+\epsilon+\frac{K'}{2(N-n)}\right)+N\alpha^2 \notag\\
			&+\sqrt{\frac{N}{2}}\alpha^2t_0\left(\sqrt{n(N-n+4)}L_1+\frac{L_2}{\sqrt{\epsilon}}+2L_3\right).
		\end{align}

		By the definition of $\mathcal{F}$ and $t_0<t$, it yields that 
		\begin{align}
			F^2(\nabla f)-\alpha f_t 
			&\le \frac{N\alpha^2}{t}+\frac{N\alpha^2}{2}\Bigg(L_1+\frac{K-\epsilon}{\alpha-1}+\frac{K'}{2(\alpha-1)(N-n)} \notag\\
			&\quad\quad\quad\quad\quad+\sqrt{\frac{2n(N-n+4)}{N}}L_1+\sqrt{\frac{2}{\epsilon N}}L_2+2\sqrt{\frac{2}{N}}L_3\Bigg).
		\end{align}
		It finishes the proof of the global gradient estimates. 
		\end{proof}

		Next, we give the proof of the Harnack inequality, i.e., Theorem \ref{thm2}.

		\begin{proof}[Proof of Theorem \ref{thm2}]
			It follows from Theorem \ref{thm1} that
			\begin{align}
				\partial_t f\le -\frac{1}{\alpha}F_2(\nabla f)+\frac{N\alpha}{t}+\frac{N\alpha}{2}\mathcal{Q}(N,n,L_1,L_2,L_3,\alpha,K, K',\epsilon),
			\end{align}
			where $\mathcal{Q}$ is defined as in \eqref{eq-Q}.

		Let $\eta(s)$ be a smooth curve connecting from $x_2$ to $x_1$ with $\eta(1)=x_1,\eta(0)=x_2$, and $F(\dot{\eta}(s))|_{\xi}$ be the length of $\dot{\eta}(s)$ at time  $\xi(s)=(1-s)t_2+t_1$, let $f(s)=\log u\left(\eta(s),\xi(s)\right)$. Then 
		\begin{align}
			f(x_1,t_1)-f(x_2,t_2)&=\int_0^1\frac{d}{ds}\left[f(\eta(s),\xi(s))\right]ds\notag\\
			&=\int_0^1(t_2-t_1)\left(\frac{df(\dot{\eta}(s))}{t_2-t_1}-\partial_t f\right)ds\notag\\
			&\le \int_0^1(t_2-t_1)\left\{\frac{F(\dot{\eta}(s))F(\nabla f)}{t_2-t_1}-\partial_t f\right\}ds\notag\\
			&\le \int_0^1(t_2-t_1)\left\{\frac{F(\dot{\eta}(s))F(\nabla f)}{t_2-t_1} +\frac{N\alpha}{\xi}+\frac{N\alpha}{2}\mathcal{Q}-\frac{1}{\alpha^2}F^2(\nabla f)\right\}ds\notag\\
			&\le  \int_0^1\frac{\alpha F^2(\dot{\eta}(s))|_{\xi}}{4(t_2-t_1)}ds+N\alpha \log \frac{t_2}{t_1}+\frac{N\alpha}{2}\mathcal{Q}(t_2-t_1).
		\end{align}

		Therefore, we arrive at 
		\begin{equation}
			u(x_1,t_1)\le u(x_2,t_2)\left(\frac{t_2}{t_1}\right)^{N\alpha} \exp\left\{ \int_0^1\frac{\alpha}{4}\frac{F^2(\dot{\eta}(s))|_{\xi}}{t_2-t_1}ds +\frac{N\alpha}{2}\mathcal{Q}(t_2-t_1)\right\}.
		\end{equation}
	\end{proof}


\section*{Acknowledgments}

The author would like to thank Professor Chen, Bin for his helpful discussion. The author is grateful to the reviewers for their careful reviews and valuable comments. 

\section*{Declarations}

\subsection*{Ethical Approval}
Ethical Approval is not applicable to this article as no human or animal studies in this study.

\subsection*{Funding} 
The author is supported partially by the NNSFC (Nos. 12001099, 12271093).

\subsection*{Data availability statement}
Data sharing is not applicable to this article as no new data were created or analyzed in this study.

\subsection*{Materials availability statement}
Materials sharing is not applicable to this article.

\subsection*{Conflict of interest/Competing interests}
All authors disclosed no relevant relationships.


{\small
	\begin{thebibliography}{99}
		
		\bibitem{Perelman:2002} G. Perelman, The entropy formula for the Ricci flow and its geometric applications, arXiv:math/0211159 [math.DG].
		
		\bibitem{Ni:2002} L. Ni, The Poisson equation and Hermitian-Einstein metrics on holomorphic vector bundles over complete noncompact Kähler manifolds, \emph{Indiana Univ. Math. J.} \textbf{51} (2002), no. 3, 679--704.156.
		
		\bibitem{Ni-Tam:2004} L. Ni and L.-F. Tam, K\"ahler-Ricci flow and the Poincaré-Lelong equation.\emph{Comm. Anal. Geom.} \textbf{12} (2004), no. 1-2, 111--141.
		
		\bibitem{Li:2024} X. Li, Entropy and heat kernel on generalized Ricci flow. \emph{J. Geom. Anal.} \textbf{34} (2024), No. 2, Paper No. 42, 23 p.
		
		\bibitem{Kristaly:2023} A. Krist\'aly, \'A. Mester, I.-I. Mezei, Sharp Morrey-Sobolev inequalities and eigenvalue problems on Riemannian-Finsler manifolds with nonnegative Ricci curvature. \emph{Commun. Contemp. Math.}  \textbf{25} (2023), No. 10, 2250063, 27 p.
		
		\bibitem{Ohta:2014} S. Ohta,  K.-T. Sturm, Bochner-Weitzenb\"ock formula and Li-Yau estimates on Finsler manifolds. \emph{Adv. Math.} \textbf{252} (2014), 429--448.
		
		\bibitem{Xia:2023} Q. Xia, Li–Yau’s Estimates on Finsler Manifolds. \emph{J. Geom. Anal.} \textbf{33} (2023), 2,  Paper No. 49, 33 pp.
		
		\bibitem{XiaC:2014} C. Xia, Local gradient estimate for harmonic functions on Finsler manifolds. \emph{Calc. Var. PDE.} \textbf{51} (2014), 3-4, 849--865.
		
		\bibitem{Bao:2007} D. Bao, On two curvature-driven problems in Riemann-Finsler geometry. \emph{Finsler Geometry in memory of Makoto Matsumoto, Adv. Stud. Pure Math.} \textbf{48}, Math. Soc., Japan, Tokyo, 2007, 19--71.
		
		\bibitem{Lakzian:2015} S. Lakzian, Differential Harnack estimates for positive solutions to heat equation under Finsler-Ricci flow, \emph{Pacific J. Math.} \textbf{278} (2015), 447462pp.
		
		\bibitem{Zeng-He:2019} F. Zeng and Q. He, Gradient estimates for a nonlinear heat equation under the Finsler-Ricci flow, \emph{Math. Slovaca} \textbf{69} (2019), No. 2, 409424.
		
		\bibitem{Zeng:2020} F. Zeng, Gradient estimates for a nonlinear heat equation under the Finsler geometric flow. \emph{J. PDE.} \textbf{33} (2020), 1, 17–-38.
		
		\bibitem{Zeng:2021} F. Zeng, Gradient estimates and Harnack inequalities for a nonlinear heat equation with the Finsler Laplacian. \emph{J. Math. Phys. Anal. Geom.} \textbf{17} (2021), 4, 521–-548.
		
		\bibitem{Cheng:2022} X. Cheng, Gradient estimates for positive solutions of heat equations under Finsler-Ricci flow, \emph{J. Math. Anal. Appl.}, \textbf{508} (2022), 2, 125897pp.
		
		\bibitem{Azami:2023} S. Azami, Differential Harnack Inequality for a Parabolic Equation Under the Finsler-Geometric Flow, \emph{Mediterr. J. Math.} (2023) 20:58.
		
		\bibitem{Y.Shen:2016} Y. Shen and Z. Shen, Introductions to modern Finsler geometry, \emph{World Scientific Publishing Company}, Singapore, 2016.
		
		\bibitem{Ohta:2009} S. Ohta, Finsler interpolation inequalities. \emph{Calc. Var. PDE.} \textbf{36} (2009), 2, 211--249.
		
		
\end{thebibliography}}
\end{document}